\documentclass[12pt]{amsart}

\usepackage{amssymb,amsmath,amsthm,enumerate,graphicx,stmaryrd}

\newcommand{\excise}[1]{}

\newcommand{\fa}{\mathfrak{a}}
\newcommand{\fb}{\mathfrak{b}}
\newcommand{\fc}{\mathfrak{c}}
\newcommand{\ox}{\mathcal{O}_X}
\newcommand{\oxp}{\mathcal{O}_{X'}}
\def\<{\langle}
\def\>{\rangle}
\def\AA{\mathcal{A}}

\def\CC{\mathbb{C}}
\def\II{\mathcal{I}}
\def\JJ{\mathcal{J}}
\def\NN{\mathbb{N}}
\def\OO{\mathcal{O}}
\def\PP{\mathcal{P}}
\def\QQ{\mathbb{Q}}
\def\RR{\mathbb{R}}
\def\ZZ{\mathbb{Z}}

\def\xx{\mathbf{x}}
\def\onto{\twoheadrightarrow}
\def\spot{{\hbox{\raisebox{.8pt}{\large\bf \hspace{-.2pt}.}}\hspace{-.5pt}}}
\def\floor#1{{\lfloor #1 \rfloor}}

\DeclareMathOperator{\link}{link}

\DeclareMathOperator{\Spec}{Spec}

\theoremstyle{plain}
\newtheorem{thm}{Theorem}[section]
\newtheorem{prop}[thm]{Proposition}
\newtheorem{lemma}[thm]{Lemma}
\newtheorem{cor}[thm]{Corollary}
\newtheorem{theorem}{Theorem}
\newtheorem{corollary}[theorem]{Corollary}
\newtheorem{proposition}[theorem]{Proposition}

\theoremstyle{definition}
\newtheorem{defn}[thm]{Definition}
\newtheorem{remark}[thm]{Remark}
\newtheorem{example}[thm]{Example}

\begin{document}

\mbox{}
\vspace{-1.1ex}
\title{Multiplier ideals of sums via cellular resolutions}
\author{Shin-Yao Jow and Ezra Miller}
\date{10 April 2007}

\begin{abstract}
Fix nonzero ideal sheaves $\fa_1,\dots,\fa_r$ and $\fb$ on a normal
$\QQ$-Gorenstein complex variety~$X$.  For any positive real numbers
$\alpha$ and~$\beta$, we construct a resolution of the multiplier
ideal $\JJ((\fa_1+\cdots+\fa_r)^\alpha \fb^\beta)$ by sheaves that are
direct sums of multiplier ideals
$\JJ(\fa_1^{\lambda_1}\cdots\fa_r^{\lambda_r}\fb^\beta)$ for various
$\lambda \in \RR^r_{\geq 0}$ satisfying $\sum_{i=1}^r \lambda_i =
\alpha$.  The resolution is cellular, in the sense that its boundary
maps are encoded by the algebraic chain complex of a regular
CW-complex.  The CW-complex is naturally expressed as a
triangulation~$\Delta$ of the simplex of nonnegative real vectors
$\lambda \in \RR^r$ with $\sum_{i=1}^r \lambda_i = \alpha$.  The
acyclicity of our resolution reduces to that of a cellular free
resolution, supported on~$\Delta$, of a related monomial ideal.  Our
resolution implies the multiplier ideal sum formula
\[
  \JJ(X,(\fa_1+\cdots+\fa_r)^\alpha \fb^\beta) =
  \sum_{\lambda_1+\cdots+\lambda_r=\alpha}
  \JJ(X,\fa_1^{\lambda_1}\cdots\fa_r^{\lambda_r}\fb^\beta),
\]
generalizing Takagi's formula for two summands \cite{Takagi}, and
recovering Howald's multiplier ideal formula for monomial
ideals~\cite{Howald} as a special case.  Our resolution also yields a
new exactness proof for the \emph{Skoda complex}
\mbox{\cite[Section~9.6.C]{Laz}}.
\end{abstract}


\maketitle

\section*{Introduction}

Let $X$ be a smooth complex algebraic variety and let $\fa\subseteq
\OO_X$ be an ideal sheaf.  Applications in algebraic geometry of the
multiplier ideal sheaves
\[
  \JJ(\fa^\alpha) = \JJ(X,\fa^\alpha) \subseteq \ox
\]
for real numbers $\alpha > 0$ have led to investigations of their
behavior with respect to natural algebraic operations.  For example,
Demailly, Ein, and Lazarsfeld \cite{DEL} proved that given two ideal
sheaves $\fa_1$ and~$\fa_2$, one has
\[
  \JJ((\fa_1\fa_2)^{\alpha})
  \subseteq
  \JJ(\fa_1^{\alpha})\JJ(\fa_2^{\alpha}).
\]
For the subtler case of sums, on the other hand, Musta\c{t}\v{a}
\cite{Mu} showed that
\[
  \JJ((\fa_1+\fa_2)^\alpha)
  \subseteq
  \sum_{0\leq t \leq \alpha}\JJ(\fa_1^{\alpha-t})\JJ(\fa_2^t),
\]
and Takagi \cite{Takagi} later refined this to
\begin{equation}\label{Takagi}
  \JJ((\fa_1+\fa_2)^\alpha)
  =
  \sum_{0\leq t \leq \alpha}\JJ(\fa_1^{\alpha-t}\fa_2^{t}),
\end{equation}
where he proved it more generally when $X$ is normal and
$\QQ$-Gorenstein.

Takagi used characteristic $p$ methods to deduce~\eqref{Takagi}, which
makes his work distinctly different from \cite{DEL} and \cite{Mu},
where the arguments proceed by geometric techniques such as log
resolutions and sheaf cohomology.  Our purpose is to show that such
geometric techniques, combined with combinatorial methods from
topology and commutative algebra, can recover Takagi's
equality~\eqref{Takagi} and generalize it.  As a consequence, we
derive a new proof (Corollary~\ref{c:3}) of Howald's formula for
multiplier ideals of monomial ideals \cite{Howald}, and demonstrate
how it can be reformulated to hold for all ideals.  In addition, we
obtain a new ``cellular'' exactness proof (Corollary~\ref{c:skoda})
for the \emph{Skoda complex} \cite[Section~9.6.C]{Laz} via the
contractibility of simplices.  Our main results center around the
following, which constitutes part of Theorem~\ref{t:JJ}.

\begin{theorem} \label{t:1}
Fix nonzero ideal sheaves $\fa_1,\dots,\fa_r,\fb$ on a normal
$\QQ$-Gorenstein complex variety~$X$ and $\alpha,\beta > 0$.  There is
a resolution \mbox{$0 \rightarrow \JJ_{r-1} \rightarrow \dots
\rightarrow \JJ_0 \rightarrow 0$} of the multiplier ideal
$\JJ((\fa_1+\cdots+\fa_r)^\alpha \fb^\beta)$ by sheaves $\JJ_i$ that
are finite direct sums of multiplier ideals of the form
$\JJ(\fa_1^{\lambda_1}\cdots\fa_r^{\lambda_r}\fb^\beta)$ for various
nonnegative $\lambda \in \RR^r$ with \mbox{$\sum_{i=1}^r \lambda_i =
\alpha$}.  Every distinct ideal sheaf of that form appears as a
summand of~$\JJ_0$.
\end{theorem}

Part of the final claim of Theorem~\ref{t:1} is that there are only
finitely many distinct multiplier ideals of the form
$\JJ(X,\fa_1^{\lambda_1}\cdots\fa_r^{\lambda_r}\fb^\beta)$ for
$\lambda_1+\cdots+\lambda_r=\alpha$.  In particular, the surjection
$\JJ_0 \onto \JJ(X,(\fa_1+\cdots+\fa_r)^\alpha \fb^\beta)$ in our
resolution implies the following (finite) summation formula; see also
Section~\ref{s:apps} for refinements.

\begin{corollary}\label{c:2}
$\displaystyle \JJ(X,(\fa_1+\cdots+\fa_r)^\alpha\fb^\beta) =
\sum_{\lambda_1+\cdots+\lambda_r=\alpha}
\JJ(X,\fa_1^{\lambda_1}\cdots\fa_r^{\lambda_r}\fb^\beta)$.
\end{corollary}

Corollary~\ref{c:2} reduces the calculation of the multiplier ideals
of arbitrary polynomial ideals to those of principal ideals.  In the
special case of a monomial ideal $\fa =
\<\xx^{\gamma_1},\ldots,\xx^{\gamma_r}\>$, generated by the monomials
in the polynomial ring $\CC[x_1,\ldots,x_d]$ with exponent vectors
$\gamma_1,\ldots,\gamma_r \in \NN^d$, the summation formula becomes
particularly explicit.  For a subset $\Gamma \subseteq \RR^d$, let
$\mathrm{conv}\,\Gamma$ denote its convex hull.  By the \emph{integer
part}\/ of a vector $\nu = (\nu_1,\ldots,\nu_d) \in \RR^d$, we mean
the vector $(\floor{\nu_1},\ldots,\floor{\nu_d}) \in \ZZ^d$ whose
entries are the greatest integers less than or equal to the
coordinates of~$\nu$.

\begin{corollary}\label{c:3}
If $\fa = \<\xx^{\gamma_1},\ldots,\xx^{\gamma_r}\>$ is a monomial
ideal in $\CC[x_1,\ldots,x_d]$, then $\JJ(\fa^\alpha)$ is generated by
the monomials in $\CC[x_1,\ldots,x_d]$ whose exponent vectors are the
integer parts of the vectors in
$\mathrm{conv}\{\alpha\cdot\gamma_1,\ldots,\alpha\cdot\gamma_r\}
\subseteq \RR^d$.
\end{corollary}
\begin{proof}
Using Corollary~\ref{c:2} with $\fa_j = \<\xx^{\gamma_j}\>$, it
suffices to note that the divisor of a single monomial has simple
normal crossings, so no log resolution is necessary.%
\end{proof}

It is easy to check that for $\alpha = 1$, the vectors in the
conclusion of Corollary~\ref{c:3} are precisely those from Howald's
result \cite{Howald}, namely the vectors $\gamma \in \NN^d$ such that
$\gamma + (1,\ldots,1)$ lies in the interior of the convex hull of all
exponents of monomials in~$\fa$.

Our approach to Theorem~\ref{t:1} is to construct a specific
resolution satisfying the hypotheses, including the part
about~$\JJ_0$.  The resolution we construct is \emph{cellular}, in a
sense generalizing the manner in which resolutions of monomial ideals
can be cellular \cite{BS}; see \cite[Chapter~4]{MS} for an
introduction.  In general, a complex in any abelian category could be
called \emph{cellular}\/ if each homological degree is a direct sum
indexed by the faces of a CW-complex, and the boundary maps are
determined in a natural way from those of the CW-complex.  An
elementary way to phrase this in the present context is as follows.
(Theorem~\ref{t:JJ} is more precise: a specific triangulation~$\Delta$
is constructed in Section~\ref{s:freeres}, and the sheaves
$\JJ_\sigma$ are specified in Remark~\ref{rk:JJ}.)

\begin{theorem} \label{t:4}
Resume the notation from Theorem~\ref{t:1}.  There is a triangulation
$\Delta$ of the simplex $\{\lambda \in \RR^r \mid \sum_{i=1}^r
\lambda_i = \alpha$ and $\lambda \geq 0\}$ such that
\[
  \JJ_i = \bigoplus_{\sigma\in\Delta_i}\JJ_\sigma
\]
can be taken to be a direct sum indexed by the set~$\Delta_i$ of
$i$-dimensional faces $\sigma \in \Delta$, and the differential
of~$\JJ_\spot$ is induced by natural maps between ideal sheaves, using
the signs from the boundary maps of~$\Delta$.  If $\lambda \in
\Delta_0$ is a vertex, then $\JJ_\lambda =
\JJ(X,\fa_1^{\lambda_1}\cdots\fa_r^{\lambda_r}\fb^\beta)$.
\end{theorem}

Comparing the final sentences of Theorems~\ref{t:1} and~\ref{t:4}, a
key point is that every possible multiplier ideal of the form
$\JJ(X,\fa_1^{\lambda_1}\cdots\fa_r^{\lambda_r}\fb^\beta)$ occurs at
some vertex $\lambda \in \Delta_0$.  Writing down what it means for
two such multiplier ideals to coincide, this stipulation provides
strong hints as to potential choices for triangulations; see
Section~\ref{s:freeres}.

The proof of exactness for the cellular resolution in
Theorem~\ref{t:4} proceeds by lifting the problem to an appropriate
log resolution $X' \to X$; see the proof of Theorem~\ref{t:JJ}.
Over~$X'$, we resolve the lifted ideal sheaf by a complex (of locally
principal ideal sheaves in~$\oxp$) that is, analytically locally at
each point of~$X'$, a cellular free complex over a polynomial ring.
This cellular free complex turns out to be a cellular free resolution
of an appropriate monomial ideal; its construction and proof of
acyclicity occupy Section~\ref{s:freeres}, particularly
Proposition~\ref{p:freeres}.  Having cellularly resolved the lifted
sheaf over~$X'$, the desired cellular resolution over~$\ox$ is
obtained by pushing forward to~$X$ and using local vanishing; again,
see the proof of Theorem~\ref{t:JJ}.  Thus Theorems~\ref{t:1}
and~\ref{t:4} constitute a certain global version of cellular free
resolutions of monomial ideals.

The acyclicity of the cellular free resolutions in
Proposition~\ref{p:freeres} reduce to a simplicial homology vanishing
statement, Corollary~\ref{c:M}, for simplicial complexes obtained by
deleting boundary faces from certain contractible
manifolds-with-boundary.  We deduce Corollary~\ref{c:M} from the
following more general statement, which is of independent interest.
Its \emph{rimmed}\/ hypothesis is satisfied by the barycentric
subdivision of any polyhedral homology-manifold-with-boundary.  In
what follows, to \emph{delete a simplex}\/~$\sigma$ from a simplicial
complex~$M$ means to remove from~$M$ every simplex
containing~$\sigma$.

\begin{proposition}\label{p:M}
Fix a simplicial complex~$M$ whose geometric realization $|M|$ is a
homology-manifold with boundary~$\partial M$.  Assume that~$M$ is
\emph{rimmed}, meaning that
\begin{equation}\label{*}
\text{if $\sigma$ is a face of~$M$, then $\sigma \cap |\partial M|$ is
a face of~$M$.}
\end{equation}
Then deleting any collection of boundary simplices from~$M$ results in
a simplicial subcomplex whose (reduced) homology is canonically
isomorphic to that of~$M$.
\end{proposition}

\subsection*{Acknowledgements}
We are grateful to Rob Lazarsfeld for valuable discussions, and to
Shunsuke Takagi for the suggestion to work with singular
varieties.  EM was supported by NSF CAREER grant DMS-0449102 and a
University of Minnesota McKnight Land-Grant Professorship.  The latter
funded his visit to the University of Michigan, which he wishes to
thank for its hospitality while this work was completed.

\section{Combinatorial topological preliminaries}\label{s:topological}

In this section we collect some elementary results from simplicial
topology.  For additional background, see \cite{munkres}, especially
\S 63 for homology-manifolds, and \cite[Section~5.1]{Ziegler} for
polyhedral complexes.  The reader interested solely in the
construction of resolutions for Theorems~\ref{t:1} and~\ref{t:4}, as
opposed to their acyclicity proofs, can proceed to
Section~\ref{s:freeres} after Example~\ref{ex:PP}.  The goal for the
remainder of this section is the proof of Proposition~\ref{p:M} and
its immediate consequence, Corollary~\ref{c:M}.

Our convention is to identify any abstract simplicial complex~$\Delta$
with the underlying topological space of any geometric
realization~$|\Delta|$.  Thus, fixing a vertex set~$V$, we view
$\Delta$ as a collection of \emph{simplices} (finite subsets of~$V$),
called \emph{faces}\/ of~$\Delta$, such that any subset of any face
of~$\Delta$ is a face of~$\Delta$.  For example, we can express the
\emph{link}\/ in~$\Delta$ of a face~$\sigma$ as the subcomplex
\[
  \link_\Delta(\sigma)=\{\tau\in \Delta \mid \tau\cup\sigma \in
  \Delta\text{ and } \tau\cap\sigma=\varnothing\}.
\]
In what follows, all links are taken inside of the ambient simplicial
complex, unless otherwise noted.  Our motivation is the following
class of simplicial complexes.

\begin{example}\label{ex:PP}
Let $\PP$ be a polyhedral complex, such as a polyhedral subdivision of
a polytope.  (This example works just as well for any regular cell
complex; see \cite{bjorner}).  The \emph{barycentric subdivision}
of~$\PP$ is the simplicial complex~$\Delta$ whose vertex set is the
set of faces of~$\PP$, and whose simplices are the chains $P_1 <
\cdots < P_\ell$ of faces $P_j \in \PP$.  Here $P < Q$ means that $P$
is a proper face of~$Q$.  The barycentric subdivision~$\Delta$ is
homeomorphic to (the space underlying)~$\PP$; one way to realize
$\Delta$ is to let each face $P_1 < \cdots < P_\ell$ be the convex
hull of the barycenters of the faces $P_1,\ldots,P_\ell$.
\end{example}

The main idea in the proof of Proposition~\ref{p:M} is to delete the
boundary simplices one by one, in a suitable order, so that at each
step the homology remains unchanged, using the following.

\begin{lemma}\label{l:filter}
Fix a simplicial complex~$\Delta$ and a chain of simplicial
subcomplexes
\[
  \Delta = \Delta_0 \supseteq \Delta_1 \supseteq \cdots \supseteq
  \Delta_r = \Gamma.
\]
If the relative chain complex of each pair $(\Delta_{i-1},\Delta_i)$
has vanishing homology for all $i = 1,\ldots,r$, then the homology
of\/~$\Gamma$ is canonically isomorphic to that of~$\Delta$.
\end{lemma}
\begin{proof}
Repeatedly apply the long exact sequence of homology.
\end{proof}

Throughout the rest of this section, we use $M$ to denote a simplicial
homology-manifold-with-boundary of dimension $\dim(M) = d$, defined as
follows.

\begin{defn}
A \emph{simplicial homology-manifold-with-boundary}\/ of dimension~$d$
is a simplicial complex whose maximal faces all have dimension~$d$,
and such that the link of each $i$-face has the homology of either a
ball or a sphere of dimension $d-i-1$.
\end{defn}

\begin{remark}\label{r:M}
If the polyhedral complex~$\PP$ in Example~\ref{ex:PP} subdivides a
homology-manifold-with-boundary, then the conditions of
Proposition~\ref{p:M}, including the rimmed condition~(\ref{*}), are
satisfied by the barycentric subdivision of~$\PP$.
\end{remark}

\begin{lemma}\label{l:link}
The link of any $i$-face in $M$ is a dimension $d-i-1$ simplicial
homology-manifold-with-boundary.
\end{lemma}
\begin{proof}
If $\sigma \subseteq M$ is a face of dimension~$i$, then $\Gamma =
\link_M(\sigma)$ is pure of dimension~\mbox{$d-i-1$} because $M$ is
pure of dimension~$d$.  The condition on the homology of the link in
$\Gamma$ of a face $\tau \in \Gamma$ holds simply because the
condition holds for \mbox{$\link_M(\sigma\cup\tau)$}, which equals
$\link_\Gamma(\tau)$ by definition.
\end{proof}

The above lemma required no special hypotheses on~$M$.  Our final
observation before the proof of Proposition~\ref{p:M} is that the
condition~(\ref{*}) of being rimmed is inherited by links of boundary
faces.

\begin{lemma}\label{l:link*}
Under the hypotheses of Proposition~\ref{p:M}, the link of any
boundary simplex of~$M$ satisfies all of the hypotheses on~$M$ in
Proposition~\ref{p:M}, including being rimmed.
\end{lemma}
\begin{proof}
Let $\sigma$ be a boundary simplex of~$M$.  In view of
Lemma~\ref{l:link}, we need only show that $\link(\sigma)$
satisfies~(\ref{*}).  This condition is a consequence of the equality
\begin{equation}\label{**}
  \partial(\link(\sigma)) = \partial M \cap \link(\sigma),
\end{equation}
because for any simplex $\tau$ in~$\link(\sigma)$, the intersection
$\tau \cap \partial(\link(\sigma))$ is forced to equal $\tau \cap
\partial M$, which is a face of~$\tau$ by~(\ref{*}).

To prove~(\ref{**}), we first show that if $\tau\in
\partial(\link(\sigma))$ then $\tau\in \partial M$; note that for this
implication, condition~(\ref{*}) is not necessary.  Replacing $\tau$
by a face of $\partial(\link(\sigma))$ containing it, if necessary, we
may assume that $\dim(\sigma\cup\tau) = d-1$; equivalently, $\tau$ is
a maximal face of $\partial(\link(\sigma))$.  In this case, there is
only one maximal face of $\link(\sigma)$ containing~$\tau$.  But the
faces of $\link(\sigma)$ containing~$\tau$ are in bijection with the
faces of~$M$ containing $\sigma\cup\tau$.  Hence there is only one
maximal face of~$M$ containing $\sigma\cup\tau$.  Therefore
$\sigma\cup\tau$ must be a boundary face of~$M$.  We conclude that
$\tau \in \partial M$, because $\tau$ is a face of $\sigma\cup\tau$.

For the reverse containment, suppose now that $\tau\in\partial M \cap
\link(\sigma)$, and let $\tau' = \tau\cup\sigma \in M$.  The
intersection $\tau'\cap \partial M$ is a face of~$M$ by
condition~(\ref{*}), and it contains all of the vertices of $\sigma$
and~$\tau$, because both $\sigma$ and~$\tau$ are boundary faces
of~$M$.  The only face of $\tau'$ containing all of the vertices of
$\sigma$ and~$\tau$ is $\tau'$ itself; hence $\tau' = \tau' \cap
\partial M$ is a boundary face of~$M$.  It follows that $\link(\tau')$
has the homology of a ball (rather than a sphere).  But $\link(\tau')
= \link_M(\sigma \cup \tau) = \link_{\link(\sigma)}(\tau)$, so $\tau
\in \partial(\link(\sigma))$.
\end{proof}

\begin{proof}[Proof of Proposition~\ref{p:M}]
Use induction on the dimension of~$M$.  For $\dim M = 1$ the statement
is an elementary claim concerning subdivided intervals, so assume that
$\dim M=d \geq 2$.  Let $S$ be the collection of boundary faces to be
deleted.  Without loss of generality, we may assume that $S$ is a
\emph{cocomplex} inside $\partial M$, which means that if $\sigma\in
S$ and $\tau\in\partial M$ such that $\tau \supseteq \sigma$, then
$\tau \in S$.

Totally order the simplices in~$S$ in such a way that all of the
simplices of maximal dimension $d-1$ come first, then those of
dimension $d-2$, and so on.  Let $\sigma_1,\sigma_2,\ldots$ be this
totally ordered sequence of simplices.  Let $M_i$ be the result of
deleting $\sigma_1,\ldots,\sigma_i$ from~$M$, with $M = M_0$.  The
desired result will follow from Lemma~\ref{l:filter}, with \mbox{$M =
\Delta$}, as soon as we show that the relative chain complexes
$\mathcal{C}_\spot(M_{i-1},M_i)$ are all acyclic.

Let $\Gamma$ be the simplicial complex obtained by deleting the
boundary from $\link_M(\sigma_i)$, which is acyclic by induction, via
Lemma~\ref{l:link*}.  We claim that $\mathcal{C}_\spot(M_{i-1},M_i)$
is (noncanonically) isomorphic to the reduced chain complex
$\tilde{\mathcal{C}}_\spot(\Gamma)$.  Indeed, recall the bijective
correspondence between the faces of $\link_M(\sigma_i)$ and the faces
of~$M$ containing~$\sigma_i$.  Under this correspondence, $\Gamma$ is
mapped bijectively to the set
\[
  \{\tau\in M \mid \tau \cap \partial M = \sigma_i\}
  =
 \{\tau \in M_{i-1} \mid \tau \supseteq \sigma_i\}.
\]
The faces of this cocomplex inside $M_{i-1}$ constitute a free basis
of $\mathcal{C}_\spot(M_{i-1},M_i)$, and this induces an isomorphism
from $\tilde{\mathcal{C}}_\spot(\Gamma)$ to
$\mathcal{C}_\spot(M_{i-1},M_i)$.
\end{proof}

For future reference, here is the special case of
Proposition~\ref{p:M} that we apply later.

\begin{cor}\label{c:M}
If $M$ is a contractible rimmed simplicial manifold-with-boundary,
then deleting any collection of boundary faces from~$M$ results in a
simplicial subcomplex with vanishing reduced~\mbox{homology}.
\end{cor}

\begin{remark}
It would be preferable to conclude in Corollary~\ref{c:M} that the
subcomplex is contractible, and more generally that the subcomplex in
Proposition~\ref{p:M} is homotopy-equivalent to~$M$, using
\cite[Lemma~10.3(i)]{bjorner}: if $\Delta_i = \Delta_{i-1} \cup
\Gamma_i$ is a union of two simplicial subcomplexes such that
$\Gamma_i$ and $\Delta_{i-1}\cap\Gamma_i$ are both contractible, then
$\Delta_i$ is homotopy-equivalent to~$\Delta_{i-1}$.  The notation
here is consistent with Lemma~\ref{l:filter}; in our case
$\Delta_{i-1}$ is the result of deleting a boundary simplex~$\sigma_i$
from~$\Delta_i$, and $\Gamma_i$ is the cone from~$\sigma_i$ over
$\Delta_{i-1} \cap \Gamma_i$ (thus $\Gamma_i$ is the \emph{closed
star} of~$\sigma_i$ in~$\Delta_i$).  The problem is that, while
$\Gamma_i$ is always contractible (it is a cone), the subcomplex
$\Delta_{i-1} \cap \Gamma_i = \link_{\Delta_i}(\sigma_i)$ need not be
simply-connected, even though it has vanishing homology.
\end{remark}

\section{A cellular free resolution}\label{s:freeres}

Let $A \in \ZZ^{n \times r}$ be an integer matrix with $n$~rows and
$r$~columns, and fix a real column vector $b \in \RR^n$ with
coordinates $b_1,\ldots,b_n$.  (When applying these results in
Section~\ref{s:log2cell}, all entries in the matrix~$A$ will be
nonnegative.)  Viewing the rows $A^1,\ldots,A^n$ as functionals
on~$\RR^r$, we get an (infinite but locally finite) affine hyperplane
arrangement
\[
  \AA = \bigcup_{\substack{z\in\ZZ^n\\1 \leq j \leq n}}
  \{\lambda\in\RR^r\mid A^j\cdot\lambda+b_j=z_j\}.
\]
The arrangement~$\AA$ induces a polyhedral subdivision of~$\RR^r$.
Fixing a nonnegative real number $\alpha \in \RR_{\geq 0}$, the
restriction of this polyhedral subdivision to the simplex
\[
  \Delta^\alpha = \{(\lambda_1,\ldots,\lambda_r) \in \RR_{\geq 0}^r
  \mid \lambda_1 + \cdots + \lambda_r = \alpha\}
\]
is a polyhedral subdivision~$\AA_\alpha$ of~$\Delta^\alpha$.  Two
points $\lambda$ and~$\mu$ lie in the same (relatively open)
\emph{cell}\/ of~$\AA_\alpha$ if and only if for each $j =
1,\ldots,n$, there is an integer~$z_j$ such that~either
\begin{itemize}
\item
$A^j\cdot\lambda+b_j = A^j\cdot\mu+b_j = z_j$, or else
\item
$A^j\cdot\lambda+b_j$ and $A^j\cdot\mu+b_j$ both lie in the
open interval $(z_j, z_j + 1)$.
\end{itemize}

For the duration of this paper, fix any polyhedral subdivision $\PP$
of~$\Delta^\alpha$ that \emph{refines}~$\AA_\alpha$, meaning that
every face of~$\PP$ is contained in a face of~$\AA_\alpha$.  As in
Example~\ref{ex:PP}, denote by $\Delta$ the barycentric subdivision
of~$\PP$, realized so that the vertices of~$\Delta$ are the
barycenters of the cells of~$\PP$.  We write $\sigma(P)$ for the
vertex of~$\Delta$ that is the barycenter of the polytope $P \in \PP$.

Construct a module $M \subseteq \CC[x_1^{\pm1},\ldots,x_n^{\pm1}]$
generated by (Laurent) monomials over the polynomial ring $S =
\CC[x_1,\ldots,x_n]$, so $M$ is a \emph{Laurent monomial module}
\cite[Definition~9.6]{MS}, as follows.  For each polytope $P \in \PP$,
set
\[
  m_P = x_1^{\floor{A^1\cdot\lambda + b_1}} \cdots
  x_n^{\floor{A^n\cdot\lambda+b_n}}\text{ for any }\lambda\in P^\circ,
\]
where $\floor{z}$ is the greatest integer less than or equal to~$z$,
and the cell $P^\circ$ is the relative interior of~$P$; the monomial
$m_P$ is well-defined because the greatest integer parts defining it
are equal for all vectors in the same cell of~$\PP$.  Now define
\[
  M = \<m_P \mid P \in \PP\>.
\]

We can label each simplex in $\Delta$ with a monomial as follows.  By
virtue of the fact that each polytope in~$\PP$ indexes a
monomial~$m_P$ in the Laurent polynomial ring, the vertices of the
simplicial complex~$\Delta$, which are the barycenters $\sigma(P)$ of
the polytopes $P \in \PP$, are labeled with monomials.  For each such
vertex~$\sigma(P)$, let us set $m_{\sigma(P)} = m_P$.  More generally,
for each simplex $\sigma \in \Delta$,
\[
  \text{if } \sigma = (P_1 < \cdots < P_\ell) \text{ then set }
  m_\sigma = m_{P_1};
\]
i.e., we label $\sigma$ by the monomial of the smallest face in the
chain.  For (Laurent) monomials~$m$ and~$m'$, we say that $m$
\emph{divides}~$m'$ if $m'/m$ is an honest monomial in~$S$.

\begin{lemma}\label{l:lcm}
If $P < Q$ are polytopes in~$\PP$, then $m_Q$ divides~$m_P$.
Consequently, if $\sigma \in \Delta$ is a simplex, then $m_\sigma$ is
the least common multiple of $\{m_\tau \mid \tau$ is a vertex
of~$\sigma\}$.
\end{lemma}
\begin{proof}
We prove only the first sentence, as the second follows easily.
{}From the perspective of the greatest integer parts used in
constructing $m_P$ and~$m_Q$, the only difference between the
barycenters of~$P$ and~$Q$ is that more of the affine functions
$\lambda \mapsto A^j\cdot\lambda+b_j$ can take integer values
on~$\sigma(P)$ than on~$\sigma(Q)$.  In particular, if
$\floor{A^j\cdot\sigma(Q)+b_j} = z_j$, then
$\floor{A^j\cdot\sigma(P)+b_j}$ equals either $z_j$ or $z_j + 1$.
This being true for all~$j$, we conclude that $m_Q$ divides~$m_P$ (in
fact, $m_P/m_Q$ is a squarefree monomial in~$S$).
\end{proof}

Lemma~\ref{l:lcm} is precisely the condition under which the monomial
labels $m_\sigma$ for $\sigma \in \Delta$ define a \emph{cellular free
complex}\/ on~$\Delta$
\cite[Section~1]{BS} (see also Chapters~4 and~9 in \cite{MS} for
background).  Briefly, this is the complex
\[
  0 \to \bigoplus_{\sigma\in\Delta_{r-1}} \<m_\sigma\> \to \cdots \to
  \bigoplus_{\sigma\in\Delta_i} \<m_\sigma\> \to \cdots \to
  \bigoplus_{\sigma\in\Delta_0} \<m_\sigma\> \to 0
\]
of free $S$-modules in which the morphism $\<m_\sigma\> \to
\<m_\tau\>$ of principal Laurent monomial modules is, for each
codimension~1 face $\tau \subset \sigma$, the natural inclusion times
$\pm 1$.  The sign is determined by arbitrary but fixed orientations
for the simplices in~$\Delta$.

\begin{prop}\label{p:freeres}
The cellular free complex supported on $\Delta$, in which $\sigma \in
\Delta$ is labeled by $m_\sigma$, is a cellular free resolution of the
Laurent monomial module~$M$.
\end{prop}
\begin{proof}
By
\cite[Proposition~1.2]{BS}, we need to show that the simplicial
subcomplex
\[
  \Delta_{\preceq m} = \{\sigma \in \Delta \mid m_\sigma \text{
  divides } m\}
\]
is acyclic for all Laurent monomials~$m$.  Let $m = x_1^{e_1} \cdots
x_n^{e_n}$.  For \mbox{$\sigma = (P_1 < \cdots < P_\ell)$}, the
monomial $m_\sigma$ divides~$m$ if and only if the barycenter $\tau =
\sigma(P_1)$ satisfies the conditions $A^j\cdot\tau + b_j < e_j + 1$
for all~$j$.

Compare the subcomplex $\Gamma_{\preceq m} \subseteq \Delta$
consisting of those points $\lambda \in |\Delta|$ such that $A^j
\cdot \lambda + b_j \leq e_j + 1$ for all~$j$.  This subcomplex is a
(compact convex) polytope, and hence $\Gamma_{\preceq m}$ is a
contractible simplicial manifold-with-boundary.  But $\Delta_{\preceq
m}$ is obtained from $\Gamma_{\preceq m}$ by deleting those
vertices~$\tau$, necessarily on the boundary of $\Gamma_{\preceq m}$,
having $A^j \cdot \tau + b_j = e_j + 1$ for some~$j$.  The desired
result is therefore a consequence of Corollary~\ref{c:M}, which
applies by Remark~\ref{r:M}.
\end{proof}

\section{Log resolution to cellular resolution}\label{s:log2cell}

For the duration of this section, fix a normal $\QQ$-Gorenstein
complex variety~$X$ and ideal sheaves $\fa_1,\ldots,\fa_r,\fb$ on~$X$.
In addition, fix a log resolution $\pi \colon X'\rightarrow X$ of
$\fb$ and $\fa_i+\fa_j$ for all $i,j\in \{1,\ldots,r\}$.  By
definition, this means that $\fb \cdot \oxp$ and $(\fa_i + \fa_j)
\cdot \oxp$ are locally principal ideal sheaves on~$X'$.  The special
case $\fa_i = \fa_j$ implies that
\begin{align*}
  \fa_i\cdot\oxp &=\oxp(-A_i) \text{ for } i \in \{1,\dots,r\}
  \text{ and }
\\
  \fb\cdot\oxp &=\oxp(-B)
\end{align*}
are the ideal sheaves associated to effective divisors $A_i$ and~$B$
on~$X'$.

We wish to calculate the multiplier ideals $\JJ((\fa_1 + \cdots +
\fa_r)^\alpha\fb^\beta)$ for nonnegative real numbers $\alpha$
and~$\beta$.  This requires a log resolution of $\fa = \fa_1 + \cdots
+ \fa_r$.  Fortunately, the morphism $\pi: X' \to X$ is one, as can be
seen by applying the following locally on~$X'$.

\begin{lemma} \label{divisor}
Let $R$ be a local integral domain.
\begin{enumerate}[1.]
\item
If $\<f_1,\ldots,f_n\>$ is a nonzero principal ideal of~$R$, then
there is an index $i \in \{1,\ldots,n\}$ such that $f_i$ divides $f_j$
for all $j \in \{1,\ldots,n\}$.
\item
If $f_1,\ldots,f_r \in R$ are elements such that $\<f_i,f_j\>$ is a
principal ideal for all $i,j \in \{1,\ldots,r\}$, then there exists a
permutation $\xi \in S_r$ such that $f_{\xi(1)} | f_{\xi(2)} | \cdots
| f_{\xi(r)}$.  In particular $\<f_1,\dots,f_r\> = \<f_{\xi(1)}\>$ is
a principal ideal.
\end{enumerate}
\end{lemma}
\begin{proof}
Suppose $\<f_1,\ldots,f_n\>$ is generated by $f \in R$.  Then each
$f_j = f_j'\cdot f$ is a multiple of~$f$.  On the other hand, $f =
\sum c_jf_j = f \cdot \sum c_j f_j'$ with $c_j \in R$.  Since $R$ is a
domain, this implies that $1 = \sum c_j f_j'$.  Since $R$ is local, we
conclude that some $f_i'$ is a unit.  This proves the first claim.
The second is an immediate consequence of the first.
\end{proof}

Similar to $\fa_i$ and $\fb$ from before, the ideal sheaf $\fa \cdot
\oxp = \oxp(-C)$ is determined by an effective divisor $C$ on~$X'$.
Let $D_1,\ldots,D_n$ be the set of prime divisors appearing in
$A_1,\ldots,A_r$, $B$, and the relative canonical divisor $K_{X'/X}$.
Thus
\begin{align*}
  A_j &= \sum_{i=1}^n a_{ij} D_i \text{ for } j \in \{1,\dots,r\},
\\
  \beta B - K_{X'/X} &= \sum_{i=1}^n b_i D_i, \text{ and }
\\
  C &= \sum_{i=1}^n \gamma_i D_i
\end{align*}
in terms of $D_1,\ldots,D_n$, where the coefficients $a_{ij}$ and
$\gamma_i$ are integers, but $b_i$ can be real numbers.  Recall that
the multiplier ideal of~$\fa$ with weighting coefficient~$\alpha$ is
\[
  \JJ(\fa^\alpha)
  =
  \pi_*\oxp(-\floor{\alpha C - K_{X'/X}})
  \subseteq
  \ox,
\]
where $\floor{\sum_{i=1}^n \delta_i D_i} = \sum_{i=1}^n
\floor{\delta_i} D_i$ is the greatest integer divisor less than or
equal to $\sum_{i=1}^n \delta_i D_i$.  More generally, for $\lambda
\in \RR_{\geq 0}^r$ and $\beta \geq 0$, set
\[
  \JJ(\fa_1^{\lambda_1}\cdots\fa_r^{\lambda_r}\fb^\beta)
  =
  \pi_*\oxp(-\floor{\lambda_1 A_1+\cdots+\lambda_r A_r+\beta B-K_{X'/X}})
  \subseteq
  \ox.
\]

How do these multiplier ideals depend on the choices of~$\lambda \in
\RR^r$ and $\beta \geq 0$?  Given the decompositions into prime
divisors $D_1,\ldots,D_n$ with integer coefficient matrix $A =
(a_{ij})$ and real vector $b = (b_i)$, we are verbatim in the
situation from Section~\ref{s:freeres}, whose notation we assume
henceforth.

\begin{lemma}\label{l:JJ}
For each polytope $P \in \PP$, there is a single ideal sheaf $\JJ_P
\subseteq \ox$ that coincides with the multiplier ideals
$\JJ(\fa_1^{\lambda_1}\cdots\fa_r^{\lambda_r}\fb^\beta)$ for all
$\lambda \in P^\circ$ in the interior of~$P$.
\end{lemma}
\begin{proof}
In fact, the greatest integer divisors $\floor{\lambda_1
A_1+\cdots+\lambda_r A_r+\beta B - K_{X'/X}}$ on~$X'$ all coincide for
$\lambda \in P^\circ$, by construction.  Therefore we can define
\[
  \II_P = \oxp(-\floor{\lambda_1 A_1 +\cdots+ \lambda_r A_r + \beta B
  - K_{X'/X}})
\]
for any $\lambda \in P^\circ$, and $\JJ_P = \pi_*\II_P$ does not
depend on $\lambda \in P^\circ$.
\end{proof}

As we did for monomials in Section~\ref{s:freeres} (before
Lemma~\ref{l:lcm}), having defined objects $\II_P$ and~$\JJ_P$ indexed
by polytopes in~$\PP$, we can define objects $\II_\sigma$
and~$\JJ_\sigma$ indexed by simplices in~$\Delta$:
\[
  \text{if } \sigma = (P_1 < \cdots < P_\ell) \text{ then set }
  \II_\sigma = \II_{P_1},
\]
and set $\JJ_\sigma = \pi_*\II_\sigma$.
\begin{remark}\label{rk:JJ}
Unwinding the definitions, we find that
\[
  \JJ_\sigma = \JJ(\fa_1^{\lambda_1}\cdots\fa_r^{\lambda_r}\fb^\beta),
\]
where $\lambda$ is the barycenter of the smallest polytope~$P_1$ in
the chain $P_1 < \cdots < P_\ell$ corresponding to the simplex $\sigma
\in \Delta$.  (Recall that $\Delta$ is the barycentric subdivision of
the polyhedral complex~$\PP$ from the beginning of
Section~\ref{s:freeres}.)
\end{remark}
The analogue of Lemma~\ref{l:lcm} holds here, as well.

\begin{lemma}\label{l:lcm'}
If $P < Q$ are polytopes in~$\PP$, then $\II_Q \supseteq \II_P$.
Consequently, if $\sigma \in \Delta$, then $\II_\sigma = \bigcap
\{\II_\tau \mid \tau$ is a vertex of~$\sigma\}$.  The same is true
with $\JJ$ in place of~$\II$.
\end{lemma}
\begin{proof}
The argument is essentially the same as the proof of
Lemma~\ref{l:lcm}.
\end{proof}

\begin{lemma}\label{l:sum}
Resume the notation involving $\fa = \fa_1+\cdots+\fa_r$ from after
Lemma~\ref{divisor}, so $\fa\cdot\oxp = \oxp(-C)$.  Then, as ideal
sheaves over~$X'$,
\[
  \oxp(-\floor{\alpha C + \beta B - K_{X'/X}}) = \sum_{P\in\PP} \II_P.
\]
\end{lemma}
\begin{proof}
Let $p \in X'$ be an arbitrary point.  By Lemma~\ref{divisor}, there
is an open neighborhood $U$ of~$p$ and a permutation $\xi \in S_r$
such that after restricting to~$U$, we have $A_{\xi(1)} \preceq \cdots
\preceq A_{\xi(r)}$, where $E \preceq F$ for divisors $E$ and~$F$
means that $F - E$ is effective.  It follows that on $U$ we have $C =
A_{\xi(1)}$, and also
\[
  \alpha A_{\xi(1)} \preceq \lambda_1 A_1 + \cdots + \lambda_r A_r
\]
whenever $\lambda_1 + \ldots + \lambda_r = \alpha$.  Hence both sides
of the desired equality are equal to $\oxp(-\floor{\alpha A_{\xi(1)} +
\beta B - K_{X'/X}})$ on~$U$.  Equality holds on~$X'$ because $p$ is
arbitrary.
\end{proof}

We have arrived at our main result.  For notational purposes, choose
an arbitrary but fixed collection of orientations for the simplices
in~$\Delta$.  There results an incidence function that assigns to each
codimension~$1$ face $\tau$ of each simplex~$\sigma$ a sign
$(-1)^{\sigma,\tau}$.

\begin{thm}\label{t:JJ}
With $\JJ_\sigma$ as in Remark~\ref{rk:JJ}, there is an exact sequence
\[
  0 \to \bigoplus_{\sigma\in\Delta_{r-1}} \JJ_\sigma \to \cdots \to
  \bigoplus_{\sigma\in\Delta_i} \JJ_\sigma \to \cdots \to
  \bigoplus_{\sigma\in\Delta_0} \JJ_\sigma \to \JJ((\fa_1 + \cdots +
  \fa_r)^\alpha\fb^\beta) \to 0
\]
of sheaves on~$X$, in which the morphism $\JJ_\sigma \to \JJ_\tau$ is
inclusion times $(-1)^{\sigma,\tau}$.  Here, $\Delta_i = \{\sigma \in
\Delta \mid \dim(\sigma) = i\}$ is the set of $i$-dimensional
simplices in~$\Delta$.
\end{thm}
\begin{proof}
We shall first verify the exactness of the complex
\[
  0 \to \bigoplus_{\sigma\in\Delta_{r-1}} \II_\sigma \to \cdots \to
  \bigoplus_{\sigma\in\Delta_i} \II_\sigma \to \cdots \to
  \bigoplus_{\sigma\in\Delta_0} \II_\sigma \to
  \sum_{\sigma\in\Delta_0} \II_\sigma \to
  0
\]
on~$X'$, where it should be noted that the sum at the far right is not
a direct sum.    It is enough to verify exactness at the stalk of an
arbitrary point $p \in X'$.  Reordering the prime divisors
$D_1,\ldots,D_n$ appearing in
$A_1,\ldots,A_r$, $B$, and $K_{X'/X}$ if necessary, we assume that
only $D_1,\ldots,D_k$ pass through~$p$.  Since $D_1,\ldots,D_n$
intersect transversally, if we pass to an analytic neighborhood or
completion at~$p$, then we can choose local coordinates
$x_1,\ldots,x_k$ so that $D_i$ is given by the vanishing of~$x_i$ for
$i = 1,\ldots,k$.  In these coordinates, $\II_\sigma$ becomes the
principal Laurent monomial module
in~$\CC[x_1^{\pm1},\ldots,x_k^{\pm1}]$ generated by the
monomial~$m_\sigma$ from before Lemma~\ref{l:lcm}, except that all of
the variables $x_{k+1},\ldots,x_n$ have been set equal to~$1$
in~$m_\sigma$.  The exactness thus follows from
Proposition~\ref{p:freeres}, given that our polyhedral complex~$\PP$
refines the subdivision induced by the linear functionals
$A^j\cdot\lambda+b_j$ for $j = 1,\ldots,k$, instead of $j =
1,\ldots,n$.  (Another reason the exactness follows from
Proposition~\ref{p:freeres} is that setting the variables equal to~$1$
is the exact operation of ``homogeneous localization''; see
\cite[Section~3.6]{Miller}, particularly Proposition~3.31.2 there.)

Pushing forward under~$\pi$ completes the proof.  Indeed,
Lemma~\ref{l:sum} implies that
\[
  \sum_{\sigma\in\Delta_0} \II_\sigma = \sum_{P\in\PP} \II_P =
  \oxp(-\floor{\alpha C + \beta B - K_{X'/X}})
\]
pushes forward to yield $\JJ((\fa_1 + \cdots +
\fa_r)^\alpha\fb^\beta)$, while local vanishing
\cite[Theorem~9.4.17(i)]{Laz} guarantees that the higher direct images
vanish on all~$\II_\sigma$.
\end{proof}

\begin{remark}
In the special case where $X=\Spec \CC[x_1,\ldots,x_n]$ and
$\fa_1,\ldots,\fa_r,\fb$ are all principal monomial ideals, all of the
$\JJ_\sigma$ are principal monomial ideals as well, so the exact
sequence in Theorem~\ref{t:JJ} becomes a cellular resolution of the
monomial ideal $\JJ((\fa_1 + \cdots + \fa_r)^\alpha\fb^\beta)$ in the
sense of \cite[Definition~4.3]{MS}.
\end{remark}

\begin{example}
Let us illustrate Theorem~\ref{t:JJ} by an example.  Take $X=\Spec
\CC[x,y]$, $r=\alpha=2$, $\fa_1=\<xy\>$, $\fa_2=\<x+y\>$, and
$\fb=\ox$.  A direct computation shows that on the doubled $1$-simplex
$\{(\lambda_1,\lambda_2)\in \RR_{\geq 0}\mid \lambda_1+\lambda_2=2\}$,
we have
\[
\JJ(\fa_1^{\lambda_1}\fa_2^{\lambda_2}) =
\begin{cases}
  \<x^2y^2\>   &\text{if } \lambda_2 = 0     \\
  \<xy\>       &\text{if } 0 < \lambda_2 < 1 \\
  \<xy(x+y)\>  &\text{if } \lambda_2 = 1     \\
  \<x+y\>      &\text{if } 1 < \lambda_2 < 2 \\
  \<(x+y)^2\>  &\text{if } \lambda_2 = 2.
\end{cases}
\]
These regions where $\JJ(\fa_1^{\lambda_1}\fa_2^{\lambda_2})$ stays
constant yield a polyhedral subdivision $\PP$ of the doubled
$1$-simplex, and we let $\Delta$ be the barycentric subdivision
of~$\PP$.  Explicitly, $\Delta$ consists of five vertices
\[
  v_1=(2,0),\ v_2=({\textstyle\frac32,\frac12}),\ v_3=(1,1),\
  v_4=({\textstyle\frac12,\ \frac32}),\ v_5=(0,2)
\]
and four edges
\[
  l_i = \text{the line segment between } v_i \text{ and }v_{i+1},
  \text{ for } i=1,2,3,4.
\]
The vertices in $\Delta$ are naturally labeled with the multiplier
ideals
\begin{align*}
  \JJ_{v_1} &= \JJ(\fa_1^{2}\fa_2^{0}) = \<x^2y^2\>,\\
  \JJ_{v_2} &= \JJ(\fa_1^{3/2}\fa_2^{1/2}) = \<xy\>,\\
  \JJ_{v_3} &= \JJ(\fa_1^1\fa_2^1) = \<xy(x+y)\>,\\
  \JJ_{v_4} &= \JJ(\fa_1^{1/2}\fa_2^{3/2}) = \<x+y\>,\\
  \JJ_{v_5} &= \JJ(\fa_1^0\fa_2^2) = \<(x+y)^2\>.
\end{align*}
The edges of~$\Delta$ are labeled by $\JJ_{l_i}$, which for
$i=1,2,3,4$ is defined to be the least common multiple of~$\JJ_{v_i}$
and~$\JJ_{v_{i+1}}$.  It is easily verified that $\JJ_{l_i}$ equals
the smaller one of~$\JJ_{v_i}$ and~$\JJ_{v_{i+1}}$.  The exact
sequence of Theorem~\ref{t:JJ} is
\[
  0\to \bigoplus_{i=1}^4 \JJ_{l_i}\to \bigoplus_{i=1}^5 \JJ_{v_i}\to
  \JJ((\fa_1+\fa_2)^2)\to 0,
\]
or more explicitly,
\[
\newcommand{\op}{\oplus}
  0
  \to
  \begin{array}{r@{\:}l}
  &\<x^2y^2\>\\\op&\<xy(x+y)\>\\\op&\<xy(x+y)\>\\\op&\<(x+y)^2\>
  \end{array}
  \to
  \begin{array}{r@{\:}l}
  &\<x^2y^2\>\\\op&\<xy\>\\\op&\<xy(x+y)\>\\\op&\<x+y\>\\\op&\<(x+y)^2\>
  \end{array}
  \to
  \mathcal{J}(\<xy,x+y\>^2)
  \to
  0.
\]
After canceling redundant terms, this is essentially a Koszul
resolution
\[
  0 \to \<xy(x+y)\> \to \<xy\> \oplus \<x+y\> \to
  \mathcal{J}(\<xy,x+y\>^2)\to 0.
\]
\end{example}

\section{Applications}\label{s:apps}

The Introduction already contains some immediate applications of
Theorem~\ref{t:JJ}, namely the Takagi-style summation formula in
Corollary~\ref{c:2} and the derivation of Howald's monomial multiplier
ideal formula in Corollary~\ref{c:3}.  In this section, we collect
further applications: a new exactness proof for the Skoda complex, and
additional summation formulas for multiplier ideals, including the
graded system case.

We begin with the summation formulas.  First, we note the following.

\begin{remark}
Theorem~\ref{t:JJ} still holds when $\fb^\beta$ is replaced by
$\fb_1^{\beta_1}\cdots\fb_s^{\beta_s}$, since the contants $\beta_k$
merely translate the affine hyperplane arrangement~$\AA$ in
Section~\ref{s:freeres}.
\end{remark}

\begin{cor} \label{summation formula}
Let $\fa_1,\dots,\fa_r, \fb_1,\dots,\fb_s \subseteq \OO_X$ be nonzero
ideal sheaves on a normal $\QQ$-Gorenstein variety~$X$ and let
$\alpha,\beta_1,\dots,\beta_s$ be positive real numbers.  Then
\[
  \JJ(X,(\fa_1+\dots+\fa_r)^\alpha
  \fb_1^{\beta_1}\cdots\fb_s^{\beta_s})
  =
  \sum_{\lambda_1+\dots+\lambda_r=\alpha}
  \JJ(X,\fa_1^{\lambda_1}\cdots\fa_r^{\lambda_r}
  \fb_1^{\beta_1}\cdots\fb_s^{\beta_s}).
\]
\end{cor}
\begin{proof}
This follows immediately from the surjectivity on the right of the
exact sequence in Theorem~\ref{t:JJ}.
\end{proof}

Our next result concerns the notion of \emph{graded system of ideals}.
For the definition of this and the multiplier ideal associated to it,
see \cite[Section~11.1.B]{Laz}.

\begin{cor}\label{graded systems}
Corollary~\ref{summation formula} still holds when the ideal sheaves
$\fa_i$ and $\fb_j$ are all replaced by graded systems of ideals.
\end{cor}
\begin{proof}
The formula for $r > 2$ can be obtained by repeatedly applying the
formula for $r=2$, so it suffices to prove this case, i.e.,
\[
  \JJ(X,(\fa_\spot+\fb_\spot)^\alpha
  \fc_{1,\spot}^{\beta_1}\cdots\fc_{s,\spot}^{\beta_s})
  =
  \sum_{0\leq t\leq\alpha}\JJ
  (X,\fa_\spot^{\alpha-t}\fb_\spot^{t}\cdot
  \fc_{1,\spot}^{\beta_1}\cdots\fc_{s,\spot}^{\beta_s}).
\]
First let us verify $\supseteq$, i.e.,
$
  \JJ(X,\fa_\spot^{\alpha-t}\fb_\spot^{t}\cdot
  \fc_{1,\spot}^{\beta_1}\cdots\fc_{s,\spot}^{\beta_s})\subseteq
  \JJ(X,(\fa_\spot+\fb_\spot)^\alpha
  \fc_{1,\spot}^{\beta_1}\cdots\fc_{s,\spot}^{\beta_s})
$
for any fixed $t\in [0,\alpha]$. By definition, the left-hand side is
equal to
$\JJ(X,\fa_p^{\frac{\alpha-t}{p}}\fb_p^{\frac{t}{p}}\cdot
\fc_{1,p}^{\frac{\beta_1}{p}}\cdots\fc_{s,p}^{\frac{\beta_s}{p}})$ for
some large $p$, and we see that
\begin{alignat*}{2}
  \JJ(X,\fa_p^{\frac{\alpha-t}{p}}\fb_p^{\frac{t}{p}}\cdot
  \fc_{1,p}^{\frac{\beta_1}{p}}\cdots\fc_{s,p}^{\frac{\beta_s}{p}})
& \subseteq \JJ(X,(\fa_p+\fb_p)^{\frac{\alpha}{p}}
  \fc_{1,p}^{\frac{\beta_1}{p}}\cdots\fc_{s,p}^{\frac{\beta_s}{p}})
&
& \ \quad\text{(by Corollary~\ref{summation formula})}
\\
& \subseteq
  \JJ(X,(\fa_\spot+\fb_\spot)_p^{\frac{\alpha}{p}}
  \fc_{1,p}^{\frac{\beta_1}{p}}\cdots\fc_{s,p}^{\frac{\beta_s}{p}})
&
&
\\
& \subseteq \JJ(X,(\fa_\spot+\fb_\spot)^\alpha
  \fc_{1,\spot}^{\beta_1}\cdots\fc_{s,\spot}^{\beta_s})
&
&\ \quad\text{(by definition).}
\end{alignat*}

To prove the reverse inclusion, first by definition there exists some
large $p$ such that
\[
  \JJ(X,(\fa_\spot+\fb_\spot)^\alpha
  \fc_{1,\spot}^{\beta_1}\cdots\fc_{s,\spot}^{\beta_s}) =
  \JJ(X,(\sum_{i=0}^p\fa_i\fb_{p-i})^{\frac{\alpha}{p}}
  \fc_{1,p}^{\frac{\beta_1}{p}}\cdots\fc_{s,p}^{\frac{\beta_s}{p}}),
\]
and by Corollary~\ref{summation formula} this right-hand side can be
expressed as
\[
  \sum_{\lambda_0+\dots+\lambda_p =
  \frac{\alpha}{p}}\JJ(X,(\prod_{i=1}^p\fa_i^{\lambda_i})
  (\prod_{i=0}^{p-1}\fb_{p-i}^{\lambda_i})\cdot
  \fc_{1,p}^{\frac{\beta_1}{p}}\cdots\fc_{s,p}^{\frac{\beta_s}{p}}).
\]
Now take a positive integer $m$ which is divisible by
$1,2,\dots,p$. Since $\fa_\spot$, $\fb_\spot$ and
$\fc_{j,\spot}$ are graded systems, we have
\[
  \fa_i^{\lambda_i}=\fa_i^{\frac{m}{i}\frac{i\lambda_i}{m}}\subseteq
  \fa_m^{\frac{i\lambda_i}{m}},\quad
  \fb_{p-i}^{\lambda_i}=\fb_{p-i}^{\frac{m}{p-i}\frac{(p-i)\lambda_i}{m}}\subseteq
  \fb_m^{\frac{(p-i)\lambda_i}{m}}, \quad
  \fc_{j,p}^{\frac{\beta_j}{p}}=\fc_{j,p}^{\frac{m}{p}\frac{\beta_j}{m}}\subseteq
  \fc_{j,m}^{\frac{\beta_j}{m}},
\]
so
\[
  (\prod_{i=1}^p\fa_i^{\lambda_i})(\prod_{i=0}^{p-1}\fb_{p-i}^{\lambda_i})\cdot
  \fc_{1,p}^{\frac{\beta_1}{p}}\cdots\fc_{s,p}^{\frac{\beta_s}{p}}\subseteq
  \fa_m^{\sum\limits_{i=1}^p
  \frac{i\lambda_i}{m}}\fb_m^{\sum\limits_{i=0}^{p-1}
  \frac{(p-i)\lambda_i}{m}}\cdot
  \fc_{1,m}^{\frac{\beta_1}{m}}\cdots\fc_{s,m}^{\frac{\beta_s}{m}},
\]
and since $\lambda_0+\dots+\lambda_p=\frac{\alpha}{p}$, if we let
$t:=\sum\limits_{i=0}^{p-1} (p-i)\lambda_i$, then $\sum\limits_{i=1}^p
i\lambda_i=\alpha-t$, hence
\begin{align*}
  \JJ(X,(\prod_{i=1}^p\fa_i^{\lambda_i})
  (\prod_{i=0}^{p-1}\fb_{p-i}^{\lambda_i})\cdot
  \fc_{1,p}^{\frac{\beta_1}{p}}\cdots\fc_{s,p}^{\frac{\beta_s}{p}})&\subseteq
  \JJ(X,\fa_m^{\frac{\alpha-t}{m}}\fb_m^{\frac{t}{m}}\cdot
  \fc_{1,m}^{\frac{\beta_1}{m}}\cdots\fc_{s,m}^{\frac{\beta_s}{m}})\\
  &\subseteq\JJ(X,\fa_\spot^{\alpha-t}\fb_\spot^{t}\cdot
  \fc_{1,\spot}^{\beta_1}\cdots\fc_{s,\spot}^{\beta_s}).\qedhere
\end{align*}
\end{proof}

For our final application, we assume that the ideal sheaves $\fa_1,
\ldots, \fa_r$ are locally principal on the $\QQ$-Gorenstein
variety~$X$.  For every nonnegative integer $\alpha < r$, the
inclusion $\fa_i \to \OO_X$ induces a natural map $\fa_i
\JJ(\fa^\alpha\fb^\beta) \to \JJ(\fa^{\alpha+1}\fb^\beta)$.  Writing
$\fa_\tau = \prod_{i\in\tau}\fa_i$ for $\tau \in \{0,1\}^r \subseteq
\RR^r$, the \emph{Skoda complex}\/ is the cellular complex
\[
  0 \to \fa_1\cdots\fa_r \JJ(\fb^\beta) \to \!\cdots\! \to
  \!\!\!\bigoplus_{|\tau| = r-\alpha}\!\!\! \fa_\tau \JJ(\fa^\alpha
  \fb^\beta) \to \!\cdots\!  \to \bigoplus_{i=1}^r \fa_i
  \JJ(\fa^{r-1}\fb^\beta) \to \JJ(\fa^r\fb^\beta) \to 0
\]
supported on a simplex whose faces are in bijection with the vectors
$\tau \in \{0,1\}^r$.

\begin{cor}\label{c:skoda}
The Skoda complex is exact.
\end{cor}

A proof for smooth $X$ appears in \cite[Section~9.6.C]{Laz}.  Here, we
provide an alternate proof, which works with no extra effort in the
$\QQ$-Gorenstein setting.

\begin{proof}
Consider the subdivisions $\AA_\alpha$ of the simplices
$\Delta^\alpha$ for $\alpha = 0,\ldots,r$, defined at the beginning of
Section~\ref{s:freeres}.  Every vector $\tau \in \{0,1\}^r \subseteq
\RR^r$ induces an embedding \mbox{$\Delta^\alpha \hookrightarrow
\Delta^r$}, where $\alpha = r-|\tau|$, by adding
$\tau$.  Choose a refinement $\PP$ of~$\AA_r$ so fine that the
subdivision induced on $\Delta^\alpha$ under every one of these
inclusions, for all~$\tau$ and all $\alpha = 0,\ldots,r$,
refines~$\AA_\alpha$.  For example, in the notation from the beginning
of Section~\ref{s:freeres}, let $\PP$ be the subdivision induced by
the affine hyperplane arrangement
\[
  \AA+\{0,1\}^r = \bigcup_{\substack{z\in\ZZ^n\\1 \leq j \leq n}}
  \bigcup_{\tau\in\{0,1\}^r} \{\lambda+\tau\in\RR^r\mid
  A^j\cdot\lambda+b_j=z_j\}
\]
obtained by translating the arrangement $\AA$ up by
every vector $\tau \in \{0,1\}^r$.

As in Theorem~\ref{t:JJ}, let $\Delta$ be the barycentric subdivision
of~$\PP$.  Recall the definition of $\JJ_\sigma =
\JJ(\fa_1^{\lambda_1}\cdots\fa_r^{\lambda_r}\fb^\beta)$ from
Remark~\ref{rk:JJ}, and write $\lambda(\sigma) = \lambda =
(\lambda_1,\ldots,\lambda_r)$.  Define a double complex
$\JJ_{\spot,\spot}$ in which
\[
  \JJ_{p,q} = \bigoplus_{\dim(\sigma)=q} \:
  \bigoplus_{\substack{|\tau|=p\\\tau\preceq\lambda(\sigma)}}
  \JJ_\sigma,
\]
where $\tau \preceq \lambda(\sigma)$ means that $\lambda(\sigma) -
\tau \in \Delta^{r-|\tau|}$ (equivalently, $\lambda(\sigma) - \tau$
has nonnegative entries).  Each horizontal complex $\JJ_{\spot,q}$ is
a direct sum, over $\sigma \in \Delta$ with $\dim(\sigma) = q$, of
complexes
$\JJ_\sigma \otimes_\CC \tilde C_\spot(\Gamma_\tau)$, where $\tau
\preceq \lambda(\sigma)$ is maximal and $\tilde C_\spot(\Gamma_\tau)$
is the reduced chain complex of a simplex~$\Gamma_\tau$ with $|\tau|$
many vertices.
The vertical differentials of $\JJ_{\spot,\spot}$ are induced by the
differentials of Theorem~\ref{t:JJ}, for each fixed~$\tau$.  In
particular, the vertical complex on the summands $\JJ_\sigma$ indexed
by a fixed $\tau \preceq \lambda = \lambda(\sigma)$ is a resolution of
$\fa_\tau \JJ(\fa^{r-|\tau|}\fb^\beta)$ by Theorem~\ref{t:JJ}, because
\[
  \JJ(\fa_1^{\lambda_1}\cdots\fa_r^{\lambda_r}\fb^\beta) = \fa_\tau
  \JJ(\fa_1^{\lambda_1-\tau_1}\cdots\fa_r^{\lambda_r-\tau_r}\fb^\beta).
\]
It follows that the vertical homology of~$\JJ_{\spot,\spot}$ is the
Skoda complex.  On the other hand, the horizontal homology
of~$\JJ_{\spot,\spot}$ is identically zero because every simplex
$\Gamma_\tau$ is contractible.  The standard spectral sequence
argument shows the desired result.
\end{proof}


\end{document}